\documentclass[12pt]{amsart}
\usepackage{amsmath, amsthm, amscd, amsfonts}

\newtheorem{theorem}{Theorem}[section]

\newtheorem{proposition}[theorem]{Proposition}
\newtheorem{corollary}[theorem]{Corollary}
\theoremstyle{definition}

\newtheorem{example}[theorem]{Example}

\theoremstyle{remark}

\numberwithin{equation}{section}
\begin{document}
\title{Approximate $C^*$-ternary ring homomorphisms}
\author{Mohammad Sal Moslehian}
\address{Department of Mathematics, Ferdowsi University, P. O. Box 1159, Mashhad 91775, Iran.}
\email{moslehian@ferdowsi.um.ac.ir} 
\begin{abstract}
In this paper, we establish the generalized Hyers--Ulam--Rassias
stability of $C^*$-ternary ring homomorphisms associated to the
Trif functional equation
\begin{eqnarray*}
d \cdot C_{d-2}^{l-2} f(\frac{x_1+\cdots +x_d}{d})+
C_{d-2}^{l-1}\sum_{j=1}^d f(x_j) = l \cdot \sum_{1\leq j_1< \cdots
< j_l\leq d} f(\frac{x_{j_1} + \cdots + x_{j_l}}{l}).
\end{eqnarray*}
\end{abstract}
\maketitle

\section {Introduction and preliminaries}
A {\it ternary ring of operators} (TRO) is a closed subspace of
the space $B({\mathcal H}, {\mathcal K})$ of bounded linear
operators between Hilbert spaces ${\mathcal H}$ and ${\mathcal
K}$ which is closed under the ternary product $[xyz] :=
xy^{\ast}z$. This concept was introduced by Hestenes \cite{HES}.
The class of TRO's includes Hilbert $C^*$-modules via the ternary
product $[xyz] := \langle x, y\rangle z$. It is remarkable that
every TRO is isometrically isomorphic to a corner $p {\mathcal A}
(1-p)$ of a $C^*$-algebra ${\mathcal A}$, where $p$ is a
projection. A closely related structure to TRO's is the so-called
$JC^*$-triple that is a norm closed subspace of $B({\mathcal H})$
being closed under the triple product $[xyz]=(xy^*z + zy^*x)/2$;
cf. \cite{HAR}. It is also true that a commutative TRO, i.e. a TRO
with the property $xy^*z=zy^*x$, is an associative $JC^*$-triple.

Following \cite{ZET} a {\it C*-ternary ring} is defined to be a
Banach space ${\mathcal A}$ with a ternary product $(x, y,
z)\mapsto [xyz]$ from ${\mathcal A}$ into ${\mathcal A}$ which is
linear in the outer variables, conjugate linear in the middle
variable, and associative in the sense that
$[xy[zts]]=[x[tzy]s]=[[xyz]ts]$, and satisfies
$\|[xyz]\|\leq\|x\|\|y\|\|z\|$ and $\|[xxx]\|=\|x\|^{3}$. For
instance, any TRO is a C*-ternary ring under the ternary product
$[xyz]= xy^*z$. A linear mapping $\varphi$ between $C^*$-ternary
rings is called a {\it homomorphism} if
$\varphi([xyz])=[\varphi(x)\varphi(y)\varphi(z)]$ for all
$x,y,z\in {\mathcal A}$.

The stability problem of functional equations originated from a
question of Ulam \cite{ULA}, posed in 1940, concerning the stability
of group homomorphisms. In the next year, Hyers \cite{HYE} gave a
partial affirmative answer to the question of Ulam in the context of
Banach spaces. In 1978, Th. M. Rassias \cite{RAS1} extended the
theorem of Hyers by considering the unbounded Cauchy difference
$\|f(x+y)-f(x)-f(y)\|\leq \varepsilon(\|x\|^ p+\|y\|^ p)$, where
$\varepsilon>0$ and $ p\in [0,1)$ are constants. The result of Th.
M. Rassias has provided a lot of influence in the development of
what we now call {\it Hyers--Ulam--Rassias stability} of functional
equations. In 1994, a generalization of Rassias' result, the
so-called generalized Hyers--Ulam--Rassias stability, was obtained
by G\u avruta \cite{GAV} by following the same approach as in
\cite{RAS1}. During the last decades several stability problems of
functional equations have been investigated in the spirit of
Hyers--Ulam--Rassias-G\u avruta. See \cite{CZE, H-I-R, JUN, RAS2,
MOS1} and references therein for more detailed information on
stability of functional equations.

As far as the author knows, \cite{BOU} is the first paper dealing
with stability of (ring) homomorphisms. Another related result is
that of Johnson \cite{JOH} in which he introduced the notion of
almost algebra $*$-homomorphism between two Banach $*$-algebras. In
fact, so many interesting results on the stability of homomorphisms
have been obtained by many mathematicians; see \cite{RAS3} for a
comprehensive account on the subject. In \cite{B-M} the stability of
homomorphisms between $J^*$-algebras associated to the Cauchy
equation $f(x+y)=f(x)+f(y)$ was investigated. Some results on
stability ternary homomorphisms may be found at \cite{A-M, M-S}.

Trif \cite{TRI} proved the generalized stability for the
so-called Trif functional equation
\begin{eqnarray*}
d \cdot C_{d-2}^{l-2} f(\frac{x_1+\cdots +x_d}{d})+
C_{d-2}^{l-1}\sum_{j=1}^d f(x_j) = l \cdot \sum_{1\leq j_1< \cdots
< j_l\leq d} f(\frac{x_{j_1} + \cdots + x_{j_l}}{l}),
\end{eqnarray*}
deriving from an inequality of Popoviciu \cite{POP} for convex
functions (here, $C^k_r$ denotes $\frac{r!}{k!(r-k)!}$). Hou and
Park \cite{P-H} applied the result of Trif to study
$*$-homomorphisms between unital $C^*$-algebras. Further, Park
investigated the stability of Poisson $JC^*$-algebra homomorphisms
associated with Trif's equation (see \cite{PAR1}.

In this paper, using some strategies from \cite{B-M, L-S, P-H, PAR1,
TRI}, we establish the generalized Hyers--Ulam--Rassias stability of
$C^*$-ternary homomorphisms associated to the Trif functional
equation. If a $C^*$-ternary ring $({\mathcal A}, [\;])$ has an
identity, i.e. an element $e$ such that $x = [xee] = [eex]$ for all
$x\in {\mathcal A}$, then it is easy to verify that $x\odot y :=
[xey]$ and $x^*:= [exe]$ make ${\mathcal A}$ into a unital
$C^*$-algebra (due to the fact that $\|x\odot x^*\odot x\| =
\|x\|^3$). Conversely, if $(A, \odot)$ is a (unital) $C^*$-algebra,
then $[xyz] := x\odot y^*\odot z$ makes ${\mathcal A}$ into a
$C^*$-ternary ring (with the unit $e$ such that $x\odot y = [xey]$)
(see \cite{MOS2}). Thus our approach may be applied to investigate
of stability of homomorphisms between unital $C^*$-algebras.

Throughout this paper, ${\mathcal A}$ and ${\mathcal B}$ denote
$C^*$-ternary rings. In addition, let $q=\frac{l(d-1)}{d-l}$ and $r
= -\frac{l}{d-l}$ for positive integers $l, d$ with $2\leq l\leq
d-1$. By an {\it approximate $C^*$-ternary ring homomorphism
associated to the Trif equation} we mean a mapping $f: {\mathcal
A}\to {\mathcal B}$ for which there exists a certain control
function $\varphi: {\mathcal A}^{d+3}\to [0, \infty)$ such that if
\begin{eqnarray*}
D_\mu f(x_1, \cdots, x_d, u, v, w)&=&\|d \cdot C_{d-2}^{l-2}
f(\frac{\mu x_1+\cdots + \mu x_d}{d} + \frac{[uvw]}{d \cdot
C_{d-2}^{l-2}}) + C_{d-2}^{l-1}
\sum_{j=1}^d \mu f(x_j) \\
&&- l \cdot \sum_{1\leq j_1< \cdots < j_l\leq d} \mu f(\frac{x_{j_1}
+ \cdots + x_{j_l}}{l}) - [f(u)f(v)f(w)]\|.
\end{eqnarray*}
then
\begin{eqnarray}\label{trifapp}
D_\mu f(x_1, \cdots, x_d, u, v, w)\leq\varphi(x_1, \cdots, x_d,
u, v, w),
\end{eqnarray}
for all scalars $\mu$ in a subset ${\mathbb E}$ of ${\mathbb C}$
and all $x_1, \cdots, x_d,  u, v, w\in {\mathcal A}$.

It is not hard to see that a function $T : X \to Y$ between
linear spaces satisfies Trif's equation if and only if there is an
additive mapping $S : X \to Y$ such that $T(x) = S(x) + T(0)$ for
all $x \in X$. In fact, $S(x) := (1/2)(T(x) - T(-x))$; see
\cite{TRI}.

\section{Main Results}

In this section, we are going to establish the generalized
Hyers--Ulam--Rassias stability of homomorphisms in $C^*$-ternary
rings associated with the Trif functional equation. We start our
work with investigating the case in which an approximate
$C^*$-ternary ring homomorphism associated to the Trif equation is
an exact homomorphism.

\begin{proposition} Let $T:{\mathcal A} \to {\mathcal B}$ be an approximate
$C^*$-ternary ring homomorphism associated to the Trif equation
with ${\mathbb E}={\mathbb C}$ and a control function $\varphi$
satisfying
\begin{eqnarray*}
\lim_{n\to\infty}q^{-n}\varphi(q^nx_1, \cdots, q^nx_d, q^n u,
q^nv, q^nw)=0,
\end{eqnarray*}
for all $x_1, \cdots, x_d,  u, v, w\in {\mathcal A}$. Suppose that
$T(qx)=qT(x)$ for all $x\in {\mathcal A}$. Then $T$ is a
$C^*$-ternary homomorphism.
\end{proposition}

\begin{proof}
$T(0)=0$, because $T(0)=qT(0)$ and $q>1$. We have
\begin{eqnarray*}
D_1 T(x_1, \cdots, x_d, 0, 0, 0)&=& q^{-n}D_1 T(q^nx_1, \cdots, q^nx_d, 0, 0, 0)\\
&\leq& q^{-n}\varphi(q^nx_1, \cdots, q^nx_d, 0, 0, 0).
\end{eqnarray*}
Taking the limit as $n\to\infty$ we conclude that $T$ satisfies
Trif's equation. Hence $T$ is additive. It follows from
\begin{eqnarray*}
D_\mu T(q^nx, \cdots, q^nx, 0, 0, 0) &=& q^n\|d \cdot
C_{d-2}^{l-2} (T(\mu x) -\mu T(x))\| \leq \varphi(q^nx, \cdots,
q^nx, 0, 0, 0),
\end{eqnarray*}
that $T$ is homogeneous.

Set $x_1=\cdots=x_d=0$ and replace $u, v, w$ by $q^nu, q^nv,
q^nw$, respectively, in (\ref{trifapp}). Since $T$ is homogeneous,
we have
\begin{eqnarray*}
\|T([uvw])-[T(u)T(v)T(w)]\|&=& q^{-3n}\|T([q^nu q^nv q^nw])-
[T(q^nu)T(q^nv)T(q^nw)]\| \\
&\leq& q^{-n}\varphi(0, \cdots, 0, q^nu, q^nv, q^nw),
\end{eqnarray*}
for all $u, v, w\in {\mathcal A}$. The right hand side tends to
zero as $n\to\infty$. Hence $T([uvw])=[T(u)T(v)T(w)]$ for all $u,
v, w\in {\mathcal A}$.
\end{proof}

\begin{theorem}\label{main} Let $f:{\mathcal A} \to {\mathcal B}$ be an approximate
$C^*$-ternary ring homomorphism associated to the Trif equation
with ${\mathbb E}={\mathbb T}$ and a control function $\varphi
:{\mathcal A}^{d+3} \to [0, \infty)$ satisfying
\begin{eqnarray}\label{phi}
\widetilde{\varphi}(x_1, \cdots, x_d, u, v,
w):=\sum_{j=0}^{\infty} q^{-j} \varphi(q^jx_1, \cdots, q^jx_d,
q^ju, q^jv, q^jw) < \infty ,
\end{eqnarray}
for all $x_1, \cdots, x_d, u, v, w\in{\mathcal A}$. If $f(0)= 0$, then there
exists a unique $C^*$-ternary ring homomorphism $T:{\mathcal A} \to {\mathcal B}$ such that
\begin{eqnarray*}
\|f(x) - T(x)\|\leq \frac{1}{l \cdot C_{d-1}^{l-1}}
\widetilde{\varphi}(qx, rx, \cdots, rx, 0, 0, 0),
\end{eqnarray*}
for all $x\in{\mathcal A}$.
\end{theorem}
\begin{proof} Set $u=v=w=0, \mu =1$ and replace
$x_1, \cdots ,x_d$ by $qx, rx, \cdots, rx$ in (\ref{trifapp}). Then
\begin{eqnarray*}
\|C_{d-2}^{l-1}f(qx)-l \cdot C_{d-1}^{l-1} f(x)\|\leq \varphi(qx,
rx, \cdots, rx, 0, 0, 0) \quad (x \in {\mathcal A}).
\end{eqnarray*}
One can use induction to show that
\begin{eqnarray}\label{approx}
&&\|q^{-n}f(q^nx)-q^{-m}f(q^mx)\|\nonumber\\
&\leq& \frac{1}{l \cdot C_{d-1}^{l-1}}\sum_{j=m}^{n-1}q^{-j}
\varphi\big(q^j(qx), q^j(rx), \cdots, q^j(rx), 0, 0, 0\big),
\end{eqnarray}
for all nonnegative integers $m<n$ and all $x \in {\mathcal A}$.
Hence the sequence $\{q^{-n}f(q^nx)\}_{n\in {\mathbb N}}$ is
Cauchy for all $x\in {\mathcal A}$. Therefore we can define the
mapping $T:{\mathcal A} \to {\mathcal B}$ by
\begin{eqnarray}\label{lim}
T(x) := \lim_{n\to\infty}\frac{1}{q^n} f(q^nx)\quad (x\in{\mathcal
A}).
\end{eqnarray}
Since
\begin{eqnarray*}
D_1T(x_1, \cdots ,x_d, 0, 0, 0) &=& \lim_{n\to\infty} q^{-n}D_1f(q^nx_1,\cdots, q^nx_d, 0, 0, 0)\\
&\leq& \lim_{n\to\infty} q^{-n}\varphi(q^nx_1,\cdots, q^nx_d, 0, 0, 0)\\
&=& 0,
\end{eqnarray*}
we conclude that $T$ satisfies the Trif equation and so it is
additive (note that (\ref{lim}) implies that $T(0)=0$). It
follows from (\ref{lim}) and (\ref{approx}) with $m=0$ that
\begin{eqnarray*}
\|f(x)- T(x)\| \leq \frac{1}{l \cdot C_{d-1}^{l-1}}
\widetilde{\varphi}(qx, rx, \cdots, rx, 0, 0, 0),
\end{eqnarray*}
for all $x\in {\mathcal A}$.

We use the strategy of \cite{TRI} to show the uniqueness of $T$. Let
$T'$ be another additive mapping fulfilling
\begin{eqnarray*}
\|f(x)- T'(x)\| \leq \frac{1}{l \cdot C_{d-1}^{l-1}}
\widetilde{\varphi}(qx, rx, \cdots, rx, 0, 0, 0),
\end{eqnarray*}
for all $x\in {\mathcal A}$. We have
\begin{eqnarray*}
\|T(x)- T'(x)\|&=&q^{-n}\|T(q^nx)-T'(q^nx)\|\\
&\leq& q^{-n}\|T(q^nx)-f(q^nx)\|+ q^{-n}\|f(q^nx)-T'(q^nx)\|\\
&\leq& \frac{2q^{-n}}{l \cdot C_{d-1}^{l-1}}\widetilde{\varphi}\big(q^n(qx), q^n(rx), \cdots, q^n(rx), 0, 0, 0\big)\\
&\leq& \frac{2}{l \cdot C_{d-1}^{l-1}}\sum_{j=n}^\infty q^{-j}
\varphi\big(q^j(qx),q^j(rx), \cdots, q^j(rx), 0, 0, 0\big),
\end{eqnarray*}
for all $x\in{\mathcal A}$. Since the right hand side tends to
zero as $n\to\infty$, we deduce that $T(x)=T'(x)$ for all
$x\in{\mathcal A}$.

Let $\mu\in{\mathbb T}^1$. Setting $x_1= \cdots = x_d = x$ and
$u=v=w=0$ in (\ref{trifapp}) we get
\begin{eqnarray*}
\| d \cdot C_{d-2}^{l-2} \big(f(\mu x) -\mu f(x)\big)\| \leq
\varphi(x, \cdots, x, 0, 0, 0),
\end{eqnarray*}
for all $x\in{\mathcal A}$. So that
\begin{eqnarray*}
q^{-n} \| d \cdot C_{d-2}^{l-2} \big(f(\mu q^n x) -\mu f(q^n
x)\big)\| \leq q^{-n} \varphi(q^nx, \cdots, q^nx, 0, 0, 0),
\end{eqnarray*}
for all $x\in{\mathcal A}$. Since the right hand side tends to zero as $n\to\infty$, we have
\begin{eqnarray*}
\lim_{n \to \infty}q^{-n}\|f(\mu q^n x) -\mu f(q^n x)\| = 0,
\end{eqnarray*}
for all $\mu\in{\mathbb T}^1$ and all $x\in{\mathcal A}$. Hence
\begin{eqnarray*}
T(\mu x) = \lim_{n\to \infty}\frac{f(q^n \mu x)}{q^n}= \lim_{n\to
\infty}\frac{\mu f(q^nx)}{q^n} = \mu T(x),
\end{eqnarray*}
for all $\mu\in{\mathbb T}^1$ and all $x\in{\mathcal A}$.

Obviously, $T(0x)=0=0T(x)$. Next, let $\lambda \in {\mathbb C}
\;\;(\lambda \neq 0)$, and let $M$ be a natural number greater than
$|\lambda|$. By an easily geometric argument, one can conclude that
there exist two numbers $\mu_1, \mu_2 \in {\mathbb T}$ such that
$2\frac{\lambda}{M}=\mu_1+\mu_2$. By the additivity of $T$ we get
$T\big(\frac{1}{2}x\big)=\frac{1}{2}T(x)$ for all $ x\in {\mathcal
A}$. Therefore
\begin{eqnarray*}
T(\lambda x)& = & T\big(\frac{M}{2}\cdot 2 \cdot
\frac{\lambda}{M}x\big)=MT\big(\frac{1}{2}\cdot 2\cdot
\frac{\lambda}{M}x\big) =\frac{M}{2}T\big(2\cdot \frac{\lambda}{M}x\big)\\
& = & \frac{M}{2}T(\mu_1x+\mu_2x)
=\frac{M}{2}\big(T(\mu_1x)+T(\mu_2x)\big)
\\ & = & \frac{M}{2}(\mu_1+\mu_2)T(x)
=\frac{M}{2}\cdot 2\cdot \frac{\lambda}{M}\\
&=&\lambda T(x),
\end{eqnarray*}
for all $x \in {\mathcal A}$, so that $T$ is a ${\mathbb
C}$-linear mapping.

Set $\mu =1$ and $x_1=\cdots=x_d=0$, and replace $u, v, w$ by $q^nu,
q^nv, q^nw$, respectively, in (\ref{trifapp}) to get
\begin{eqnarray*}
\frac{1}{q^{3n}}\big\|d \cdot C_{d-2}^{l-2} f\big(\frac{q^{3n}}{d \cdot C_{d-2}^{l-2}}[uvw]\big)-\big[f(q^nu)f(q^nv)f(q^nw)\big]\big\|\\
\leq q^{-3n}\varphi(0, \cdots,  0, q^nu, q^nv, q^nw),
\end{eqnarray*}
for all $u, v, w\in {\mathcal A}$. Then by applying the continuity of the ternary product
$(x,y,z)\mapsto [xyz]$ we deduce
\begin{eqnarray*}
T([uvw])&=& d \cdot C_{d-2}^{l-2} T\big(\frac{1}{d \cdot C_{d-2}^{l-2}}[uvw]\big)\\
&=&\lim_{n\to\infty}\frac{d \cdot C_{d-2}^{l-2}}{q^{3n}} f\big(\frac{q^{3n}}{d \cdot C_{d-2}^{l-2}}[uvw]\big)\\
&=&\lim_{n\to\infty}\big[\frac{f(q^nu)}{q^n}\frac{f(q^nv)}{q^n}\frac{f(q^nw)}{q^n}\big]\\
&=& [T(u)T(v)T(w)],
\end{eqnarray*}
for all $u, v, w\in {\mathcal A}$. Thus $T$ is a $C^*$-ternary
homomorphism.
\end{proof}

\begin{example} Let $S:{\mathcal A} \to {\mathcal A}$ be a (bounded) $C^*$-ternary homomorphism,
and let $f:{\mathcal A} \to {\mathcal A}$ be defined by
$$f(x)=\left \{\begin{array}{cc}S(x) \;\;\;\;\;\; \|x\|<1\\ 0
\;\;\;\;\;\;\;\;\;\;\;\; \|x\|\geq 1 \end{array} \right .$$ and
$$\varphi(x_1, \cdots, x_d, u, v, w) := \delta, $$
where $\delta := d \cdot C_{d-2}^{l-2} + d \cdot C_{d-2}^{l-1} + l
\cdot C_d^l + 1$. Then
\begin{eqnarray*}
\widetilde{\varphi}(x_1, \cdots, x_d, u, v, w)&=&\sum_{n=0}^\infty q^{-n}\cdot \delta = \frac{\delta q}{q-1},
\end{eqnarray*}
and
\begin{eqnarray*}
D_\mu f(x_1, \cdots, x_d, u, v, w)\leq \varphi (x_1, \cdots, x_d, u, v, w),
\end{eqnarray*}
for all $\mu\in {\mathbb T}^1$ and  all $x_1, \cdots, x_d,  u, v,
w\in {\mathcal A}$. Note also that $f$ is not linear. It follows
from Theorem \ref{main} that there is a unique $C^*$-ternary ring
homomorphism $T: {\mathcal A} \to {\mathcal A}$ such that
\begin{eqnarray*}
\|f(x)-T(x)\|\leq \frac{1}{l \cdot C_{d-1}^{l-1}}\,
\widetilde{\varphi}(qx, rx, \cdots, rx, 0, 0, 0) \qquad (x\in
{\mathcal A}).
\end{eqnarray*}
Further, $T(0)=\lim_{n\to\infty}\frac{f(0)}{q^n}=0$ and for $x\neq
0$ we have
\begin{eqnarray*}
T(x)=\lim_{n\to\infty}\frac{f(q^nx)}{q^n}
=\lim_{n\to\infty}\frac{0}{q^n}=0,
\end{eqnarray*}
since for sufficiently large $n, \|q^nx\|\geq 1$. Thus $T$ is
identically zero.
\end{example}

\begin{corollary}
Let $f:{\mathcal A} \to {\mathcal B}$ be a mapping with $f(0)= 0$
and there exist constants $\varepsilon \geq 0$ and $p\in[0, 1)$
such that
\begin{eqnarray*}
D_\mu f(x_1, \cdots, x_d, u, v, w)\leq \varepsilon (\sum_{j=1}^d
\|x_j\|^p  +  \|u\|^p + \|v\|^p + \|w\|^p),
\end{eqnarray*}
for all $\mu\in{\mathbb T}^1$ and all $x_1, \cdots, x_d, u, v,
w\in{\mathcal A}$. Then there exists a unique $C^*$-ternary ring
homomorphism $T:{\mathcal A} \to {\mathcal B}$ such that
\begin{eqnarray*}
\|f(x) - T(x)\|\leq \frac{q^{1-p}(q^p+(d-1)r^p)\varepsilon }{l
\cdot C_{d-1}^{l-1}(q^{1-p}-1)}\|x\|^p,
\end{eqnarray*}
for all $x\in{\mathcal A}$.
\end{corollary}

\begin{proof} Define $\varphi(x_1, \cdots, x_d, u, v, w) = \varepsilon
(\sum_{j=1}^d \|x_j\|^p + \|u\|^p + \|v\|^p + \|w\|^p)$, and apply
Theorem 2.2.
\end{proof}

The following corollary can be applied in the case that our ternary
algebra is linearly generated by its `idempotents', i.e. elements
$u$ with $u^3 = u$.
\begin{proposition}
Let ${\mathcal A}$ be linearly spanned by a set $S\subseteq
{\mathcal A}$ and let $f:{\mathcal A} \to {\mathcal B}$ be a mapping
satisfying $f(q^{2n}[s_1s_2z]) = [f(q^ns_1)f(q^ns_2)f(z)]$ for all
sufficiently large positive integers $n$, and all $s_1,s_2\in S,
z\in{\mathcal A}$. Suppose that there exists a control function
$\varphi :{\mathcal A}^{d} \to [0, \infty)$ satisfying
\begin{eqnarray*}
\widetilde{\varphi}(x_1, \cdots, x_d):=\sum_{j=0}^{\infty} q^{-j}
\varphi(q^jx_1, \cdots, q^jx_d) < \infty \quad (x_1, \cdots, x_d
\in{\mathcal A}).
\end{eqnarray*}
If $f(0)=0$ and
\begin{eqnarray*}
\|d \cdot C_{d-2}^{l-2} f(\frac{\mu x_1+\cdots + \mu x_d}{d}) +
C_{d-2}^{l-1}
\sum_{j=1}^d \mu f(x_j)\\
- l \cdot \sum_{1\leq j_1< \cdots < j_l\leq d} \mu f(\frac{x_{j_1} +
\cdots + x_{j_l}}{l})\|\leq \varphi(x_1, \cdots, x_d),
\end{eqnarray*}
for all $\mu \in{\mathbb T}^1$ and all $x_1, \cdots, x_d \in
{\mathcal A}$, then there exists a unique $C^*$-ternary ring
homomorphism $T:{\mathcal A} \to {\mathcal B}$ such that
\begin{eqnarray*}
\|f(x) - T(x)\|\leq \frac{1}{l \cdot C_{d-1}^{l-1}}\,
\widetilde{\varphi}(qx, rx, \cdots, rx),
\end{eqnarray*}
for all $x\in{\mathcal A}$.
\end{proposition}
\begin{proof}
Applying the same argument as in the proof of Theorem 2.2, there
exists a unique linear mapping $T:{\mathcal A} \to {\mathcal B}$
given by
\begin{eqnarray*}
T(x) := \lim_{n\to\infty}\frac{1}{q^n} f(q^nx) \quad
(x\in{\mathcal A})
\end{eqnarray*}
such that
\begin{eqnarray*}
\|f(x) - T(x)\|\leq \frac{1}{l \cdot C_{d-1}^{l-1}}
\widetilde{\varphi}(qx, rx, \cdots, rx),
\end{eqnarray*}
for all $x\in{\mathcal A}$. We have
\begin{eqnarray*}
T([s_1s_2z]) &=& \lim_{n\to\infty}\frac{1}{q^{2n}} f([(q^ns_1)(q^ns_2)z])\\
&=& \lim_{n\to\infty}\big[\frac{f(q^ns_1)}{q^n}\frac{f(q^ns_2)}{q^n}f(z)\big]\\
&=& [T(s_1)T(s_2)f(z)].
\end{eqnarray*}
By the linearity of $T$ we have $T([xyz]) = [T(x)T(y)f(z)]$ for all $x, y, z\in {\mathcal A}$.
Therefore $q^nT([xyz])= T([xy(q^nz)]) = [T(x)T(y)f(q^nz)]$, and so
\begin{eqnarray*}
T[xyz])=
\lim_{n\to\infty}\frac{1}{q^n}[T(x)T(y)f(q^nz)]=\big[T(x)T(y)\lim_{n\to\infty}\frac{f(q^nz)}{q^n}\big
]= [T(x)T(y)T(z)],
\end{eqnarray*}
for all $x,y,z\in{\mathcal A}$.
\end{proof}

\begin{theorem} Suppose that $f:{\mathcal A} \to {\mathcal B}$ is an approximate
$C^*$-ternary ring homomorphism associated to the Trif equation with
${\mathbb E}=\{1, {\bf i}\}$ and a control function $\varphi:
A^{d+3}\to [0, \infty)$ fulfilling (\ref{phi}). If $f(0)=0$ and for
each fixed $x\in {\mathcal A}$ the mapping $t\mapsto f(tx)$ is
continuous on ${\mathbb R}$, then there exists a unique
$C^*$-ternary homomorphism $T:{\mathcal A} \to {\mathcal B}$ such
that
\begin{eqnarray*}
\|f(x)-T(x)\|\leq \widetilde{\varphi}(qx, rx, \cdots, rx, 0, 0, 0),
\end{eqnarray*}
for all $x\in{\mathcal A}$.
\end{theorem}
\begin{proof} Put $u=v=w=0$ and $\mu=1$ in (\ref{trifapp}). Using the same argument as in the proof of Theorem
\ref{main} we deduce that there exists a unique additive mapping
$T:{\mathcal A} \to {\mathcal B}$ given by
\begin{eqnarray*}
T(x)=\lim_{n\to\infty}\frac{f(q^nx)}{q^n} \quad (x\in {\mathcal
A}).
\end{eqnarray*}
By the same reasoning as in the proof of the main theorem of
\cite{RAS1}, the mapping $T$ is ${\mathbb R}$-linear.

Putting $x_1= \cdots = x_d = x$, $\mu={\bf i}$ and $u=v=w=0$ in (\ref{trifapp}) we get
\begin{eqnarray*}
\|d \cdot C_{d-2}^{l-2} (f({\bf i} x) -{\bf i} f(x))\| \leq
\varphi(x, \cdots, x, 0, 0, 0) \quad (x\in {\mathcal A}).
\end{eqnarray*}
Hence
\begin{eqnarray*}
q^{-n}\|f(q^n{\bf i}x)-{\bf i}f(q^nx)\|\leq q^{-n}\varphi(q^nx,
\cdots, q^nx, 0, 0, 0) \quad (x\in {\mathcal A}).
\end{eqnarray*}
The right hand side tends to zero as $n\to\infty$, hence
\begin{eqnarray*}
T({\bf i}x)=\lim_{n\to\infty}\frac{f(q^n{\bf
i}x)}{q^n}=\lim_{n\to\infty}\frac{{\bf i}f(q^nx)}{q^n}={\bf i}T(x)
\quad (x\in {\mathcal A}).
\end{eqnarray*}
For every $\lambda\in {\mathbb C}$ we can write
$\lambda=\alpha_1+{\bf i}\alpha_2$ in which
$\alpha_1,\alpha_2\in{\mathbb R}$. Therefore
\begin{eqnarray*}
T(\lambda x)&=&T(\alpha_1x+{\bf i}\alpha_2x)=\alpha_1T(x)+\alpha_2T({\bf i}x)\\
&=&\alpha_1T(X)+{\bf i}\alpha_2T(x)=(\alpha_1+{\bf i}\alpha_2)T(x)\\
&=&\lambda T(x),
\end{eqnarray*}
for all $x\in {\mathcal A}$. Thus $T$ is ${\mathbb C}$-linear.
\end{proof}

\end{document}